\documentclass{article}
\usepackage{amsmath}
\usepackage{amsthm}
\usepackage{amssymb}

\usepackage{amsfonts}

\newcommand{\Z}{\mathbb{Z}}
\newcommand{\st}{\;\big | \;} 

\providecommand{\abs}[1]{\lvert #1 \rvert}

\DeclareMathOperator{\Des}{\mathrm{Des}}
\DeclareMathOperator{\des}{\mathrm{des}}
\DeclareMathOperator{\del}{\mathrm{del}}
\DeclareMathOperator{\inv}{\mathrm{inv}}
\DeclareMathOperator{\ltrm}{\mathrm{ltrm}}

\DeclareMathOperator{\maj}{\mathrm{maj}}
\DeclareMathOperator{\nrmaj}{\mathrm{nrmaj}}
\DeclareMathOperator{\Neg}{\mathrm{Neg}}
\DeclareMathOperator{\rmaj}{\mathrm{rmaj}}
\DeclareMathOperator{\Abs}{\mathrm{abs}}
\DeclareMathOperator{\dmaj}{\mathrm{dmaj}}
\DeclareMathOperator{\drmaj}{\mathrm{drmaj}}

\newtheorem{thm}{Theorem}[section]
\newtheorem{prop}[thm]{Proposition}
\newtheorem{lem}[thm]{Lemma}
\newtheorem{cor}[thm]{Corollary}
\theoremstyle{definition}
\newtheorem{defn}[thm]{Definition}
\newtheorem{rem}[thm]{Remark}
\newtheorem{exmp}[thm]{Example}

\begin{document}

\title{MacMahon-type Identities for Signed Even Permutations}

\author{Dan Bernstein\\
    Department of Mathematics\\
    The Weizmann Institute of Science\\
    Rehovot 76100, Israel\\
    \texttt{dan.bernstein@weizmann.ac.il} }

\date{May 18, 2004}

\maketitle

\begin{abstract}
MacMahon's classic theorem states that the {\it length} and {\it
major index} statistics are equidistributed on the symmetric group
$S_n$. By defining natural analogues or generalizations of those
statistics, similar equidistribution results have been obtained
for the alternating group $A_n$ by Regev and Roichman, for the
hyperoctahedral group $B_n$ by Adin, Brenti and Roichman, and for
the group of even-signed permutations $D_n$ by Biagioli. We prove
analogues of MacMahon's equidistribution theorem for the group of
signed even permutations and for its subgroup of even-signed even
permutations.
\end{abstract}

\section{Introduction}

A classic theorem by MacMahon \cite{Mac} states that two {\it
permutation statistics}, namely the {\it length} (or {\it
inversion number}) and the {\it major index}, are equidistributed
on the symmetric group $S_n$. Many refinements and generalizations
of this theorem are known today (see \cite{RR} for a brief
review). In \cite{RR}, Regev and Roichman gave an analogue of
MacMahon's theorem for the alternating group $A_n \subseteq S_n$,
and in \cite{ABR}, Adin, Brenti and Roichman gave an analogue for
the hyperoctahedral group $B_n = C_2 \wr S_n$. Both results
involve natural generalizations of the $S_n$ statistics having the
equidistribution property.

Our main result here (Proposition \ref{PR:1}) is an analogue of
MacMahon's equidistribution theorem for {\it the group of signed
even permutations} $L_n = C_2 \wr A_n \subseteq B_n$. Namely, we
define two statistics on $L_n$, the {\it $L$-length} and the {\it
negative alternating reverse major index}, and show that they have
the same generating function, hence they are equidistributed. Our
Main Lemma (Lemma \ref{LE:1}) shows that every element of $L_n$
has a unique decomposition into a descent-free factor and a
signless even factor.

In \cite{B}, Biagioli proved an analogue of MacMahon's theorem for
{\it the group of even-signed permutations} $D_n$ (signed
permutations with an even number of sign changes). Using our main
result, we prove an analogue for {\it the group of even-signed
even permutations} $(L\cap D)_n = L_n \cap D_n$ (see Proposition
\ref{PR:2}).

The rest of this paper is organized as follows: Section
\ref{SEC:prelim} contains notations and definitions of the groups
and statistics we use. Also in that section, we define canonical
presentations in $B_n$ and in $L_n$. In Section
\ref{SEC:signedeven} we prove the equidistribution identity for
$L_n$, and finally in Section \ref{SEC:evensigned} we prove the
equidistribution identity for $(L \cap D)_n$.

\section{Preliminaries}\label{SEC:prelim}
\subsection{Notation}
For an integer $a\ge 0$ we let $[a]=\{1,2,\dots,a\}$ (where
$[0]=\emptyset$).

We let $C_m$ be the cyclic group of order $m$.

Let $S_n$ be the symmetric group on $1,\dots,n$ and let $A_n
\subset S_n$ denote the alternating group.

We denote by $B_n$ the wreath product $C_2 \wr S_n$. We think of
$\sigma \in B_n$ as a bijection of $\{\pm 1, \pm2, \dots, \pm n\}$
onto itself such that $\sigma(-i)=-\sigma(i)$ for all $i\in [n]$.

If $\sigma\in B_n$ then we write $\sigma=[\sigma_1, \dots,
\sigma_n]$ to mean that $\sigma(i)=\sigma_i$ for $i\in[n]$. Call
this the {\it window} notation of $\sigma$.

We denote by $L_n$ the wreath product $C_2 \wr A_n$, the subgroup
of $B_n$ of index 2 containing the signed even permutations (which
is not to be confused with the group of even-signed permutations
mentioned in Section \ref{SEC:evensigned}).

\subsection{$S_n$}
\subsubsection{The Coxeter System of $S_n$}
$S_n$ is a Coxeter group of type $A$. The Coxeter generators are
the adjacent transpositions $S=\{s_i\}_{i=1}^{n-1}$ where
$s_i=(i,i+1)$. The defining relations are the Moore-Coxeter
relations:
\begin{align*}
(s_i s_{i+1})^3 = 1 &\quad (1 \le i < n),\\
(s_i s_j)^2 = 1 &\quad (|i-j|>1),\\
s_i^2 = 1 &\quad (1\le i < n).
\end{align*}

\subsubsection{The $S$ Canonical Presentation}
For each $1 \le j \le n-1$ define
\[
    R_j^S = \{ 1, s_j, s_j s_{j-1}, \dots, s_j s_{j-1}\dots s_1 \}
\]
and note that $R_1^S,\dots,R_{n-1}^S \subseteq S_n$.
\begin{thm}\label{TH:SCanRep}
Let $w \in S_n$. Then there exist unique elements $w_j \in R_j^S$,
$1 \le j \le n-1$, such that $w = w_1 \dots w_{n-1}$. Thus, the
presentation $w=w_1\dots w_{n-1}$ is unique.
\end{thm}
\begin{defn}
Call $w=w_1\dots w_{n-1}$ in the above theorem {\it the $S$
canonical presentation of $w \in S_n$}.
\end{defn}

\subsection{$B_n$}
\subsubsection{The Coxeter System of $B_n$}

$B_n$ is a Coxeter group of type $B$, generated by
$s_1,\dots,s_{n-1}$ together with an exceptional generator
$s_0=[-1,2,3,\dots,n]$, whose action is as follows:
\begin{align*}
    [\sigma_1,\sigma_2,\dots,\sigma_n] s_0 &= [-\sigma_1,\sigma_2,\dots,\sigma_n] \\
    s_0 [\sigma_1,\dots,\pm 1,\dots,\sigma_n] &= [\sigma_1,\dots,\mp 1,\dots,\sigma_n].
\end{align*}
The additional relations are: $s_0^2 = 1$, $(s_0 s_1)^4 = 1$, and
$s_0 s_i = s_i s_0$ for all $1<i<n$.

\subsubsection{The $B$ Canonical Presentation}
For each $0 \le j \le n-1$ define
\[\begin{split}
    R_j^B = \{& 1, s_j, s_j s_{j-1}, \dots, s_j s_{j-1}\dots s_1,s_j s_{j-1}\dots s_1 s_0,\\
                & s_j s_{j-1}\dots s_1 s_0 s_1,\dots, s_j s_{j-1}\dots s_1 s_0 s_1 \dots s_j \}
\end{split}\] and note that $R_0^B,\dots,R_{n-1}^B \subseteq B_n$.
\begin{thm}\label{TH:BCanRep}
Let $\sigma \in B_n$. Then there exist unique elements $\sigma_j
\in R_j^B$, $0 \le j \le n-1$, such that $\sigma = \sigma_0 \dots
\sigma_{n-1}$. Moreover, written explicitly $\sigma_0 \dots
\sigma_{n-1}=s_{i_1}s_{i_2}\dots s_{i_r}$ is a reduced expression
for $\sigma$, that is $r$ is the minimum length of an expression
of $\sigma$ as a product of elements in $S$.
\end{thm}
\begin{defn}
Call $\sigma=\sigma_0\dots \sigma_{n-1}$ in the above theorem {\it the $B$
canonical presentation of $\sigma \in B_n$}.
\end{defn}

To prove Theorem \ref{TH:BCanRep} we use the following property of
Coxeter groups.

\begin{prop}[See \cite{GP}, Proposition 2.1.1]\label{PR:coxeterReps}
Let $(W,S)$ be a Coxeter system, and let $J \subseteq S$. Define
$
    W_J = \langle J \rangle
$
and
$
    X_J = \{ w \in W \st \ell(sw) > \ell(w) \; \forall s\in J \}
$.
\begin{enumerate}
\item
    For each $w \in W$ there exist unique $v \in W_J$ and $x \in
    X_J$ such that $w=vx$. Moreover, $\ell(w) = \ell(v)+\ell(x)$.
\item
    For any $w \in W$, $w \in X_J$ iff $w$ is the unique element
    of minimal length in $W_J w$.
\end{enumerate}
In particular, $X_J$ is a complete set of right coset
representatives of $W_J$ in $W$.
\end{prop}
The set $X_J$ is called the set of {\it distinguished right coset
representatives} of $W_J$ in $W$.

\begin{proof}[Proof of Theorem \ref{TH:BCanRep}]
First we show that there exists a reduced expression of the above
form for $\sigma$. According to the proposition above, it suffices
to show that for every $0<j \le n$, $R_{j-1}^B$ is the set of
distinguished right coset representatives of $B_{j-1}=\langle s_0,
\dots, s_{j-2} \rangle$ in $B_j$.

First we show that $R_{j-1}^B$ consists of right coset
representatives. Let $\sigma \in B_j$, $r=\abs{\sigma^{-1}(j)}$,
and $\sigma' = [\sigma_1 \sigma_2 \dots \sigma_{r-1}
\sigma_{r+1}\dots \sigma_{j} j] \in B_{j-1}$. If $\sigma(r)=j$,
then
\[
    \sigma' s_{j-1} s_{j-2} \dots s_{r} = [\sigma_1 \sigma_2 \dots \sigma_{r-1} j \sigma_{r+1}\dots
    \sigma_{j}]=\sigma;
\]
otherwise $\sigma(r)=-j$, so
\[\begin{split}
    \sigma' s_{j-1} s_{j-2} \dots s_1 s_0 s_1 \dots s_{r-1}\\
    = [\sigma_1 \sigma_2 \dots \sigma_{r-1} (-j) \sigma_{r+1}\dots
    \sigma_{j}]=\sigma.
\end{split}\]
In particular, since $s_{j-1} s_{j-2} \dots s_{r}$ and $s_{j-1}
s_{j-2} \dots s_1 s_0 s_1 \dots s_{r-1}$ are in $R_{j-1}^B$, this
shows that $\sigma \in \sigma' R_{j-1}^B$. It follows that
$B_{j-1}R_{j-1}^B = B_j$. Since $\abs{R_{j-1}^B}= 2j =
\frac{\abs{B_j}}{\abs{B_{j-1}}}$, $R_{j-1}^B$ contains exactly one
representative of each coset.

Next we note that $\{1, s_{j-1}, \dots, s_{j-1} \dots s_1\}$ are
reduced expressions in $S_j$ and that $\{1, s_0, s_0 s_1, \dots,
s_0 \dots s_{j-1}\}$ are reduced representatives of (some of) the
cosets of the parabolic subgroup $S_j$ of $B_j$, therefore
$R_{j-1}^B$ consists of reduced expressions.

Finally, the uniqueness of the canonical presentation follows from
a counting argument:
\[
    \prod_{j=0}^{n-1}\abs{R_j^B} = \prod_{j=0}^{n-1} 2(j+1) = 2^n n! = \abs{B_n}. \qedhere
\]
\end{proof}

\begin{rem}
For $\sigma \in S_n$, the $B$ canonical presentation of $\sigma$
coincides with its $S$ canonical presentation.
\end{rem}

\begin{exmp}\label{EX:BCanRep}
Let $\sigma = [5, -1, 2, -3, 4]$, then $\sigma_4 = s_4 s_3 s_2
s_1$; $\sigma \sigma_4^{-1} = [-1,2,-3,4,5]$, therefore $\sigma_3
= 1$ and $\sigma_2 = s_2 s_1 s_0 s_1 s_2$; and finally $\sigma
\sigma_4^{-1}\sigma_3^{-1}\sigma_2^{-1}=[-1,2,3,4,5]$ so $\sigma_1
= 1$ and $\sigma_0 = s_0$. Thus $\sigma = \sigma_0 \sigma_1
\sigma_2 \sigma_3 \sigma_4 = (s_0)(1)(s_2 s_1 s_0 s_1 s_2)(1)(s_4
s_3 s_2 s_1)$.
\end{exmp}

\subsection{$A_{n+1}$}

\subsubsection{A Generating Set for $A_{n+1}$}

We let
\[
    a_i = s_1 s_{i+1} \quad (1 \le i \le n-1).
\]

The set $A = \{a_i\}_{i=1}^{n-1}$ generates $A_{n+1}$, and the
generators satisfy the relations
\begin{align*}
(a_i a_j)^2 = 1 &\quad (|i-j|>1),\\
(a_i a_{i+1})^3 = 1 &\quad (1 \le i < n-1),\\
a_i^2 = 1   &\quad (1 < i \le n-1),\\
a_1^3 = 1&.
\end{align*}

Note that $(A_{n+1},A)$ is not a Coxeter system (in fact,
$A_{n+1}$ is not a Coxeter group) as $a_1^2 \neq 1$. This set of
generators was introduced by Mitsuhashi\cite{M}.

\subsubsection{The $A$ Canonical Presentation}
For each $1 \le j \le n-1$ define
\[
    R_j^A = \{ 1,\; a_j,\; a_j a_{j-1},\; \dots,\; a_j\dots a_2,\; a_j\dots
    a_2 a_1,\; a_j \dots a_2 a_1^{-1} \}
\]
and note that $R_1^A,\dots,R_{n-1}^A \subseteq A_{n+1}$.

\begin{thm}
Let $v \in A_{n+1}$. Then there exist unique elements $v_j \in
R_j^A$, $1 \le j \le n-1$, such that $v = v_1 \dots v_{n-1}$, and
this presentation is unique.
\end{thm}
\begin{defn}
Call $v=v_1\dots v_{n-1}$ in the above theorem {\it the $A$
canonical presentation of $v \in A_{n+1}$}.
\end{defn}

\subsection{$L_{n+1} = C_2 \wr A_{n+1}$}

\subsubsection{Characterization in Terms of the $B$ Canonical Presentation}
Define the group homomorphism $\Abs:C_2 \wr S_n \to S_n$ by $
((\epsilon_1,\dots,\epsilon_n),\sigma) \mapsto \sigma$, or
equivalently, in terms of our representation of elements of $C_s
\wr S_n$ as bijections of $\{\pm 1, \dots, \pm n\}$ onto itself, $
\Abs(\sigma)(i) = \abs{\sigma(i)} $.

From this formulation one sees immediately that for any $\sigma
\in B_n$, $\Abs( \sigma s_0 ) = \Abs( \sigma )$. Thus if $\sigma =
s_{i_1}\dots s_{i_k}$, then deleting all occurrences of $s_0$ from
$s_{i_1}\dots s_{i_k}$ what remains is an expression for
$\Abs(\sigma)$. Since by definition $\Abs(L_{n+1}) = A_{n+1}$, we
have
\begin{prop}
\[
    L_{n+1} = \left\{ \sigma \in B_{n+1} \st \sigma = s_{i_1} \dots
s_{i_k}, \#\{j \st i_j \neq 0\} \text{ is even} \right\}.
\]
\end{prop}

\subsubsection{A Generating Set for $L_{n+1}$}

$L_{n+1}$ is generated by $a_1,\dots,a_{n-1}$ together with the generator $a_0
= s_0 = [-1, 2, 3, \dots,n, n+1]$. The additional relations are $a_0^2 = 1$,
$(a_0 a_1)^6 = (a_0 a_1^{-1})^6 = 1$, and $(a_0 a_i)^4 = 1$ for all $1 < i \le
n-1$.

\subsubsection{The $L$ Canonical Presentation}
Let $
    R_0^L = \{ 1, a_0, a_1 a_0 a_1^{-1}, a_0 a_1 a_0 a_1^{-1} \}
$ and for each $1 \le j \le n-1$ define
\[\begin{split}
    R_j^L = & R_j^A \cup \{ a_j a_{j-1} \dots a_2 a_1^{-1} a_0,\; a_j a_{j-1} \dots a_2 a_1^{-1} a_0 a_1^{-1}\} \\
            &\cup \{ a_j a_{j-1} \dots a_2 a_1^{-1} a_0 a_1,\;\dots,\; a_j a_{j-1}\dots a_2 a_1^{-1} a_0 a_1 a_2 \dots a_j
            \} .
\end{split}\]
For example,
\[
    R_2^L = \{ 1,\; a_2,\; a_2 a_1,\; a_2 a_1^{-1},\; a_2 a_1^{-1} a_0,\;
        a_2 a_1^{-1} a_0 a_1^{-1},\; a_2 a_1^{-1} a_0 a_1,\;
        a_2 a_1^{-1} a_0 a_1 a_2 \} .
\]
Note that $R_0^L,\dots,R_{n-1}^L \subseteq L_{n+1}$.

\begin{thm}\label{TH:LCanRep}
Let $\pi \in L_{n+1}$. Then there exist unique elements $\pi_j \in R_j^L$,
$0 \le j \le n-1$, such that $\pi = \pi_0 \dots \pi_{n-1}$, and this presentation is unique.
\end{thm}
A proof is given below.
\begin{defn}
Call $\pi=\pi_0\dots \pi_{n-1}$ in the above theorem {\it the $L$
canonical presentation of $\pi \in L_{n+1}$}.
\end{defn}

The following recursive {\bf L-Procedure} is a way to calculate
the $L$ canonical presentation:

First note that $R_0^L = L_2$ so $R_0^L$ gives the canonical
presentations of all $\pi \in L_2$.

For $n>1$, let $\pi \in L_{n+1}$, $\abs{\pi(r)}=n+1$.

If $\pi(r)=n+1$, `pull $n+1$ to its place on the right' by
\[\begin{split}
    [\dots,n+1,\dots] a_{r-1} a_r \dots a_{n-1} = [\dots\dots,n+1] \quad
    \text{ if $r>2$ },\\
    [k ,n+1,\dots] a_1^{-1} a_2 \dots a_{n-1} = [\dots\dots,n+1] \quad
    \text{ if $r=2$ },\\
    (*) \quad\quad [n+1,\dots] a_1 a_2 \dots a_{n-1} = [\dots\dots,n+1] \quad
    \text{ if $r=1$ };
\end{split}\]
and if $\pi(r)=-(n+1)$, `correct the sign' by
\[\begin{split}
     [\dots,-(n+1),\dots] a_{r-2} \dots a_1^{-1} a_0 = [n+1,\dots] \quad \text{ if
     $r>3$     }, \\
     [\ell,k,-(n+1),\dots] a_1^{-1} a_0 = [n+1,\dots] \quad \text{ if
     $r=3$     }, \\
     [k, -(n+1), \dots] a_1 a_0 = [n+1,\dots] \quad \text{ if $r=2$
     }, \\
     [-(n+1),\dots] a_0 = [n+1,\dots] \quad \text{ if
     $r=1$     },
\end{split}\]
and then `pull to the right' using $(*)$.

 This gives $\pi_{n-1} \in R_{n-1}^L$ and $\pi \pi_{n-1}^{-1}
\in L_n$. Therefore by induction $\pi = \pi_0 \dots \pi_{n-2}
\pi_{n-1}$ with $\pi_j \in R_j^L$ for all $0\le j \le n-1$.

For example, let $\pi = [3, 5, -4, 2, -1]$, then $\pi_3 = a_3 a_2
a_1$; $\pi \pi_3^{-1} = [-4,3,2,-1,5]$, therefore $\pi_2 = a_2
a_1^{-1} a_0$; next $\pi \pi_3^{-1}\pi_2^{-1} = [2,3,-1,4,5]$ so
$\pi_1 = a_1$; and finally $\pi
\pi_3^{-1}\pi_2^{-1}\pi_1^{-1}=[-1,2,3,4,5]$ so $\pi_0 = a_0$.
Thus
\[
\pi = \pi_0 \pi_1 \pi_2 \pi_3 = (a_0)(a_1)(a_2 a_1^{-1} a_0)(a_3
a_2 a_1) .
\]

\begin{proof}[Proof of Theorem \ref{TH:LCanRep}]
The L-Procedure proves the existence of such a presentation, and
the uniqueness follows by a counting argument:
\[
    \prod_{j=0}^{n-1}\abs{R_j^L} = \prod_{j=0}^{n-1} 2(j+2) = 2^n (n+1)! = 2^{n+1}\abs{A_{n+1}} = \abs{L_{n+1}}. \qedhere
\]
\end{proof}

\begin{rem}
For $\pi \in A_{n+1}$, the $L$ canonical presentation of $\pi$
coincides with its $A$ canonical presentation.
\end{rem}
\begin{rem}
The canonical presentation of $\pi \in L_{n+1}$ is not necessarily
a reduced expression. For example, the canonical presentation of
$\pi=[-3,1,-2]\in L_3$ is $\pi=(a_1 a_0 a_1^{-1}) (a_1^{-1} a_0)$
which is not reduced ($\pi = a_1 a_0 a_1 a_0 $).
\end{rem}

\subsection{$B_n$ and $L_{n+1}$ Statistics}

\begin{defn}
Let $w=[w_1,w_2,\dots,w_n]$ be a word on $\Z$. Define the {\it
inversion number} of $w$ as $\inv(w)=\#\{1\le i<j \le n \st
w_i>w_j\}$.
\end{defn}

For example, $\inv([5,-1,2,-3,4]) = 6$.

\begin{defn}
1. Let $\sigma \in B_n$, then $j\ge 2$ is a {\it l.t.r.min}
(left-to-right minimum) of $\sigma$ if $\sigma(i)>\sigma(j)$ for
all $1 \le i<j$.

2. Define $\del_B(\sigma) = \#\ltrm(\sigma) = \#\{2\le j \le n \st
\text{$j$ is a l.t.r.min of $\sigma$} \}$.
\end{defn}

For example, the left-to-right minima of $\sigma=[5,-1,2,-3,4]$ are $\{2,4\}$ so
$\del_B(\sigma)=2$.

\begin{rem}
$\del_S(w)$ was defined in \cite{RR} for $w\in S_n$. According to
Proposition 7.2 in \cite{RR}, if $w\in S_n$ then
$\del_S(w)=\del_B(w)$.
\end{rem}

\begin{defn}
Let $\sigma \in B_n$. Define
\[
    \Neg(\sigma) = \{i\in[n] \st \sigma(i)<0 \}
\]
\end{defn}

\begin{lem}\label{LE:Neg}
\[
    \Neg(\sigma^{-1}) = \{ \abs{ \sigma(i) } \st i\in [n],\; \sigma(i)<0 \}
\]
\end{lem}
\begin{proof}
For $i \in [n],$
\[
\begin{split}
    \sigma^{-1}( \abs{ \sigma(i) } ) = \begin{cases}
            i   & \text{ if $\sigma(i)>0$},     \\
            -i  & \text{ if $\sigma(i)<0$} .
        \end{cases}
\end{split}
\]
\end{proof}

\begin{rem}\label{RE:SRightMultNeg}
If $v\in S_n$ and $\sigma \in B_n$ then
\begin{align*}
    \Neg(v\sigma)
        &= \{ i\in[n] \st v(\sigma(i))<0 \} \\
        &= \{ i\in[n] \st \sigma(i)<0 \} \\
        &= \Neg(\sigma).
\end{align*}
\end{rem}

\begin{defn}
Let $\sigma \in B_n$. Define the {\it $B$-length} of $\sigma$ the
usual way, i.e. $\ell_B(\sigma)$ is the length of a $\sigma$ with
respect to the Coxeter generators of $B_n$.
\end{defn}

For example,
\[
\ell_B([5,-1,2,-3,4]) = \ell_B(s_0 s_2 s_1 s_0 s_1 s_2 s_4 s_3 s_2
s_1) = 10
\]
(see Example \ref{EX:BCanRep}).

\begin{lem}
Let $\sigma \in B_n$. Then
\begin{equation}\label{EQ:ellB}
    \ell_B(\sigma) = \inv(\sigma) + \sum_{i \in \Neg(\sigma^{-1})} i   .
\end{equation}
\end{lem}

\begin{proof}
Let $\sigma=\sigma_0 \dots \sigma_{n-1}$ be the canonical
presentation and let $\sigma' = \sigma_0 \dots \sigma_{n-2} = [b_1
b_2 \dots b_{n-1} n]$. If $\sigma_{n-1}=s_{n-1} s_{n-2} \dots
s_{r}$, $r>0$, then $\sigma = [b_1 b_2 \dots b_{r-1} n b_{r} \dots
b_{n-1}]$, so $\Neg(\sigma^{-1})=\Neg(\sigma'^{-1})$ by Lemma
\ref{LE:Neg}, and $\inv(\sigma) = \inv(\sigma')+n-r$; and if
$\sigma_{n-1}=s_{n-1}\dots s_0 \dots s_r$ then $\sigma = [b_1
\dots b_r (-n) b_{r+1} \dots b_{n-1}]$, so $\Neg(\sigma^{-1}) =
\Neg(\sigma'^{-1})\cup \{n\}$ and $\inv(\sigma) =
\inv(\sigma')+r$. In both cases
\[
\big(\inv(\sigma)+ \sum_{i\in \Neg(\sigma^{-1})} i \big) -
\big(\inv(\sigma')+ \sum_{i\in \Neg(\sigma'^{-1})} i \big) =
\ell_B(\sigma_{n-1}) .
\]
By induction on $n$ the lemma is proved.
\end{proof}

In \cite{RR}, the {\it $A$-length} of $w \in A_{n}$, $\ell_A(w)$
was defined as the length of $w$'s $A$ canonical presentation, and
it was shown to have the following property.
\begin{prop}[See \cite{RR}, Proposition 4.4]\label{PR:RR4.4}
Let $w \in A_{n}$, then
\[
    \ell_A(w) = \ell_S(w)-\del_S(w) .
\]
\end{prop}

This serves as motivation for the following definition.

\begin{defn}
Let $\sigma \in B_n$. Define the {\it $L$-length} of $\sigma$ as
\begin{equation}\label{EQ:ellL}
\ell_L(\sigma) = \ell_B(\sigma)-\del_B(\sigma)
    = \inv(\sigma)-\del_B(\sigma)+\sum_{i \in \Neg(\sigma^{-1})} i .
\end{equation}
\end{defn}

\begin{rem}
1. $\ell_L$ is {\it not} a length function with respect to any set
of generators, that is for every set of generators of $L_n$, there
exists $\pi \in L_n$ such that $\ell_L(\pi)$ is in not the length
of a reduced expression for $\pi$ using those generators. For
example, in $L_3$ we have $
     \ell_L([3,1,2]) = \ell_L([-1,2,3]) = 1
$ but $
    \ell_L([3,1,2][-1,2,3]) = \ell_L([-3,1,2]) = 3      $.

2. If $w \in A_n$ then, according to Proposition \ref{PR:RR4.4}
and the above remarks, $\ell_A(w)=\ell_L(w)$.
\end{rem}

\begin{defn}
1. The {\it S-descent set} of $\sigma \in B_n$ is defined by
\[
\Des_S(\sigma) = \{ 1\le i \le n-1 \st \ell_B(\sigma s_i) <
\ell_B(\sigma) \}.
\]

2. Define the {\it major index} of $\sigma \in B_n$ by
\[
    \maj_B(\sigma) = \sum_{i \in \Des_S(\sigma) } i .
\]

3. Define the {\it reverse major index} of $\sigma \in B_n$ by
\[
    \rmaj_{B_n}(\sigma) = \sum_{i \in \Des_S(\sigma)} (n-i) .
\]

\end{defn}

For example, if $\sigma = [5,-1,2,-3,4]$ then $\Des_S(\sigma) =
\{1,3\}$, $\maj_B(\sigma) = 4$ and $\rmaj_{B_5}(\sigma) = 6$.

\begin{rem}\label{REM:DesB}
$\Des_S(\sigma) = \{ 1 \le i \le n-1 \st \sigma(i) >
    \sigma(i+1) \}$. Indeed, by Remark \ref{RE:SRightMultNeg} and the definition of $\inv$, for $1 \le i \le n-1$
\[
\begin{split}
\ell_B(\sigma s_i) - \ell_B(\sigma)
    &= \big( \inv(\sigma s_i) + \sum_{i \in \Neg((\sigma s_i)^{-1})} i \big)
        - \big( \inv(\sigma) + \sum_{i \in \Neg(\sigma^{-1})} i \big)   \\
    &= \inv(\sigma s_i) - \inv(\sigma) \\
    &= \begin{cases}
        +1  &   \text{if $\sigma(i)<\sigma(i+1)$}, \\
        -1  &   \text{if $\sigma(i)>\sigma(i+1)$}   .
    \end{cases}
\end{split}
\]
\end{rem}
\begin{rem}
It is natural to define the $B$-descent set of $\sigma \in B_n$ as
\[
\Des_B(\sigma) = \{ 0 \le i \le n-1 \st\ell_B(\sigma s_i) <
\ell_B(\sigma) \},
\]
and it also holds that
\[
\Des_B(\sigma) = \{ 0 \le i \le n-1 \st
\sigma(i) >
    \sigma(i+1) \}
\]
where $\sigma(0):=0$. Moreover, in the definition of $\maj_B$, one
can replace $\Des_S$ with $\Des_B$. However, using the $S$-descent
set makes the definitions of $\rmaj_B$ and the $L$ analogues
(below) appear more natural.
\end{rem}

$\maj_B$ and $\rmaj_{B_n}$ are equidistributed on $B_n$, as the
following lemma shows.

\begin{lem}\label{LE:phi}
There exists an involution $\phi$ of $B_n$ satisfying the
conditions
\[
    \maj_B(\sigma) = \rmaj_{B_n}(\phi(\sigma))
\]
and
\begin{equation}\label{EQ:NegUnderPhi}
    \Neg(\sigma^{-1}) = \Neg((\phi(\sigma))^{-1})   .
\end{equation}
\end{lem}
\begin{proof}
Given $\sigma = [\sigma_1,\dots,\sigma_n] \in B_n$,
$\sigma_{i_1}<\sigma_{i_2}<\dots<\sigma_{i_n}$, let $\rho_\sigma$
be the order-reversing permutation on $\{ \sigma_1,\dots,\sigma_n
\}$, that is $\rho_\sigma(\sigma_{i_k}) = \sigma_{i_{n+1-k}}$, and
define
\[
    \phi(\sigma) = [\rho_\sigma(\sigma_n),\rho_\sigma(\sigma_{n-1}),\dots,\rho_\sigma(\sigma_1)] .
\]
Since $\rho_\sigma$ is a permutation, the letters in the window
notation of $\phi(\sigma)$ are again $\sigma_1,\dots,\sigma_n$, so
$\rho_{\phi(\sigma)}=\rho_\sigma$. Thus
\[
\begin{split}
    \phi^2(\sigma) & =
        [\rho_{\phi(\sigma)}(\rho_\sigma(\sigma_1)),\dots,\rho_{\phi(\sigma)}(\rho_\sigma(\sigma_n))]
        \\
        & = [\rho_\sigma^2(\sigma_1),\dots,\rho_\sigma^2(\sigma_n)] \\
        & = \sigma
\end{split}
\]
and by Lemma \ref{LE:Neg}, $\Neg(\sigma^{-1}) =
\Neg(\phi(\sigma)^{-1})$.

Finally,
\[
\begin{split}
    i \in \Des_S(\phi(\sigma))
        & \iff \phi(\sigma)(i) > \phi(\sigma)(i+1) \\
        & \iff  \rho_\sigma(\sigma_{n+1-i}) >
        \rho_\sigma(\sigma_{n-i}) \\
        & \iff \sigma_{n+1-i} < \sigma_{n-i} \\
        & \iff n-i \in \Des_S(\sigma)       ,
\end{split}
\]
So
\[
    \rmaj_{B_n}(\phi(\sigma)) = \sum_{i \in \Des_S(\phi(\sigma))}
    n-i = \sum_{i \in \Des_S(\sigma)}i = \maj_B(\sigma) .
\]

\end{proof}

\begin{exmp}
Let $\sigma = [5,-1,2,-3,4]$. To compute $\phi(\sigma)$, we
reverse $\sigma$ to get $[4,-3,2,-1,5]$, then apply the
order-reversing permutation on $\{-3,-1,2,4,5\}$ to get
$\phi(\sigma) = [-1,5,2,4,-3]$. Indeed we have
 $\maj_B(\sigma) = 4 = \rmaj_{B_5}(\phi(\sigma))$ and $\Neg(\sigma^{-1})=\{1,3\} =
 \Neg(\phi(\sigma)^{-1})$.
\end{exmp}

\begin{defn}
1. The {\it A-descent set} of $\pi \in L_{n+1}$ is defined by
\[
    \Des_A(\pi) = \{ 1 \le i \le n-1 \st \ell_L(\pi a_i) \le \ell_L(\pi) \},
\]
and the {\it A-descent number} of $\pi \in L_{n+1}$ is defined by
$
    \des_A(\pi) = |\Des_A{\pi}|$.

2. Define the {\it alternating reverse major index} of $\pi \in
L_{n+1}$ by
\[
    \rmaj_{L_{n+1}}(\pi) = \sum_{i \in \Des_A(\pi)} (n-i) .
\]

3. Define the {\it negative alternating reverse major index} of
$\pi \in L_{n+1}$ by
\[
    \nrmaj_{L_{n+1}}(\pi) = \rmaj_{L_{n+1}}(\pi) + \sum_{i\in \Neg(\pi^{-1})}i .
\]

\end{defn}

For example, if $\pi=[5,-1,2,-3,4]$ then $\Des_A(\pi) = \{1,2\}$,
$\rmaj_{L_5}(\pi)=5$, and $\nrmaj_{L_5}(\pi)=5+1+3=9$.

\begin{rem}
1. For $w\in A_{n+1}$, the above definitions agree with the
definitions in  \cite{RR}.

2. In general, $
    \Des_A(\pi) \neq \{ 1 \le i \le n-1 \st \pi(i)>\pi(i+1) \}$.
\end{rem}

\section{Equidistribution on $L_{n+1}$}\label{SEC:signedeven}

The following is our main result.
\begin{prop}\label{PR:1}
For every $B \subseteq [n+1]$
\begin{eqnarray*}
    \sum_{\{ \pi \in
    L_{n+1} \st \Neg(\pi^{-1})\subseteq B\} }q^{\nrmaj_{L_{n+1}}(\pi)} = \sum_{\{ \pi \in
    L_{n+1} \st \Neg(\pi^{-1})\subseteq B\} }q^{\ell_L(\pi)} \\
    = \prod_{i \in
    B}(1+q^i)\prod_{i=1}^{n-1}(1+q+\dots+q^{i-1}+2q^i) .
\end{eqnarray*}
\end{prop}

By the inclusion-exclusion principle we have:

\begin{cor}\label{COR:inex}
For every $B \subseteq [n+1]$
\begin{eqnarray*}
    \sum_{\{ \pi \in
    L_{n+1} \st \Neg(\pi^{-1}) = B\} }q^{\nrmaj_{L_{n+1}}(\pi)} = \sum_{\{ \pi \in
    L_{n+1} \st \Neg(\pi^{-1}) = B\} }q^{\ell_L(\pi)}   .
\end{eqnarray*}
\end{cor}

Note that the case $B=\emptyset$ of Proposition \ref{PR:1} is just
the case $t=1$ of the following theorem.
\begin{thm}[See \cite{RR}, Theorem 6.1(2)]\label{TH:RR6.1}
\begin{multline*}
    \sum_{w \in A_{n+1}} q^{\ell_A(w)} t^{\del_A(w)} = \sum_{w \in
    A_{n+1}} q^{\rmaj_{A_{n+1}}(w)} t^{\del_A(w)} \\
    = (1+2qt)(1+q+2q^2 t)\cdots(1+q+\dots+q^{n-2}+2q^{n-1}t).
\end{multline*}
\end{thm}

The proof of Proposition \ref{PR:1} uses the decomposition of
\[
\{ \pi \in L_{n+1} \st \Neg(\pi^{-1}) \subseteq B\}
\]
into left cosets of $A_{n+1}$, and a set of distinguished coset
representatives.

\begin{lem}\label{LE:s1LeftMultLength}
Let $\rho \in S_{n+1}$. Then
\[
\ell_L(\rho)=\ell_L(s_1 \rho).
\]
\end{lem}
\begin{proof}
\[
\inv(s_1 \rho) = \begin{cases}
    \inv(\rho)+1    & \text{if $\rho^{-1}(1)<\rho^{-1}(2)$;}\\
    \inv(\rho)-1    & \text{if $\rho^{-1}(1)>\rho^{-1}(2)$}
\end{cases}
\]
and
\[
\del_B(s_1 \rho) = \begin{cases}
    \del_B(\rho)+1  & \text{if $\rho^{-1}(1)<\rho^{-1}(2)$;}\\
    \del_B(\rho)-1  & \text{if $\rho^{-1}(1)>\rho^{-1}(2)$,}
\end{cases}
\]
therefore
\begin{equation*}
    \ell_L(\rho)=\inv(\rho)-\del_B(\rho)=\inv(s_1 \rho)-\del_B(s_1 \rho)=\ell_L(s_1 \rho).
\end{equation*}
\end{proof}

\begin{lem}[Main Lemma]\label{LE:1}
Let $\pi \in L_{n+1}$. Then there exists a unique $\sigma \in
L_{n+1}$ such that $u=\sigma^{-1}\pi \in A_{n+1}$ and
$\des_A(\sigma)=0$. Moreover, $
    \Des_A(u)=\Des_A(\pi)$, $\inv(u)-\del_S(u) = \inv(\pi)-\del_B(\pi)$, and
$
    \Neg(\pi^{-1})=\Neg(\sigma^{-1})$.

\end{lem}
\begin{proof}
Let $\sigma'\in B_{n+1}$ be the increasing word with the letters
of $\pi$. Clearly $
    \inv(\sigma')=\del_B(\sigma')=0
$ so by \eqref{EQ:ellL}, $
    \ell_L(\sigma')=\sum_{i\in \Neg(\sigma'^{-1})} i$.

For every $v \in S_{n+1}$ and $i,j\in [n+1]$,
\[
    v(i)<v(j) \iff (\sigma' v)(i) < (\sigma' v)(j) ,
\]
thus
\begin{equation}\label{EQ:SRightMultInv}
\begin{aligned}
    \inv(\sigma' v)
        &= \inv(v)
\end{aligned}
\end{equation}
and
\begin{equation}\begin{aligned}\label{EQ:SRightMultDel}
    \del_B(\sigma' v)
        &=\del_S(v) .
\end{aligned}\end{equation}
By Remark \ref{RE:SRightMultNeg}, $    \Neg((\sigma'v)^{-1}) =
\Neg(v^{-1}\sigma'^{-1})
        = \Neg(\sigma'^{-1})$.
Therefore for every $v\in S_{n+1}$,
\begin{equation}\label{EQ:addLengths}
    \begin{split}
    \ell_L(\sigma' v)
        &= \inv(\sigma' v)+
            \sum_{i \in \Neg((\sigma'v)^{-1})}i \quad - \del_B(\sigma'v)\\
        &= \sum_{i \in \Neg(\sigma'^{-1})}i \quad + \inv(v) - \del_B(v)\\
        &= \ell_L(\sigma')+\ell_L(v).
    \end{split}
\end{equation}

There are two possible cases to consider:

Case 1: $\sigma' \in L_{n+1}$. Let $\sigma=\sigma'$ and let
$u=\sigma'^{-1}\pi$.

Using \eqref{EQ:addLengths} we have for $1 \le i \le n-1$,
\begin{align*}
    \ell_L(\sigma a_i)
        &= \ell_L(\sigma' a_i) \\
        &= \ell_L(\sigma') +\ell_L(a_i)\\
        &> \ell_L(\sigma') \\
        &= \ell_L(\sigma)
\end{align*}
and
\begin{align*}
    \ell_L(\pi) - \ell_L(\pi a_i)
        &= \ell_L(\sigma u) - \ell_L(\sigma (u a_i))   \\
        &= \ell_L(\sigma' u) - \ell_L(\sigma' (u a_i))   \\
        &= \ell_L(\sigma')+\ell_L(u) - \ell_L(\sigma')-\ell_L(u a_i) \\
        &= \ell_L(u)-\ell_L(u a_i).
\end{align*}
Therefore $\des_A(\sigma)=0$ and $\Des_A(u)=\Des_A(\pi)$ as
desired. From \eqref{EQ:SRightMultInv} and
\eqref{EQ:SRightMultDel} we also get that
\begin{align*}
    \inv(\pi)-\del_B(\pi)
        &= \inv(\sigma u)-\del_B(\sigma u) \\
        &= \inv(\sigma' u)-\del_B(\sigma' u) \\
        &= \inv(u)-\del_S(u).
\end{align*}

Case 2: $\sigma' s_1 \in L_{n+1}$. Let $\sigma = \sigma' s_1$ and
let $u=s_1\sigma'^{-1} \pi$.

Using \eqref{EQ:addLengths} we have for $1 \le i \le n-1$,
\begin{align*}
    \ell_L(\sigma a_i)
        &= \ell_L(\sigma' s_{i+1}) \\
        &= \ell_L(\sigma') +\ell_L(s_{i+1}) \\
        &> \ell_L(\sigma')  \\
        &= \ell_L(\sigma') + \ell_L(s_1) &\text{($\ell_L(s_1)=0$)} \\
        &= \ell_L(\sigma' s_1)  \\
        &= \ell_L(\sigma)
\end{align*}
and, using also Lemma \ref{LE:s1LeftMultLength},
\begin{align*}
    \ell_L(\pi) - \ell_L(\pi a_i)
        &= \ell_L(\sigma' s_1 u) - \ell_L(\sigma' (s_1 u a_i)) \\
        &= \ell_L(\sigma')+\ell_L(s_1 u) - \ell_L(\sigma')-\ell_L(s_1 u a_i) \\
        &= \ell_L(s_1 u)-\ell_L(s_1 (u a_i)) \\
        &= \ell_L(u)-\ell_L(u a_i).
\end{align*}
Therefore $\des_A(\sigma)=0$ and $\Des_A(u)=\Des_A(\pi)$ as
desired. From \eqref{EQ:SRightMultInv} and
\eqref{EQ:SRightMultDel} and Lemma \ref{LE:s1LeftMultLength},
\begin{align*}
    \inv(\pi)-\del_B(\pi)
        &= \inv(\sigma' s_1 u)-\del_B(\sigma' s_1 u) \\
        &= \inv(s_1 u)-\del_S(s_1 u) \\
        &= \inv(u)-\del_S(u).
\end{align*}

In both cases, the fact that $\Neg(\pi^{-1})=\Neg(\sigma^{-1})$
follows by Remark \ref{RE:SRightMultNeg} from the fact that
$\pi^{-1}=u^{-1}\sigma^{-1}$ and $u\in A_{n+1}$.

To see that $\sigma$ is unique, suppose $\tilde \sigma \in
L_{n+1}$ satisfies $\des_A(\tilde \sigma)=0$ and $\tilde u =
\tilde \sigma^{-1}\pi \in A_{n+1}$. Then $
    0
        = \des_A ( \tilde \sigma )
        = \des_A ( \sigma u \tilde u^{-1} )
$
(since $\tilde \sigma = \sigma u \tilde u^{-1}$), so for $1 \le i
\le n-1$,
\begin{align*}
    0
        &\le \ell_L(\sigma u \tilde u^{-1} a_i) -\ell_L(\sigma u \tilde
            u^{-1}) \\
        &= \ell_L(\sigma)+\ell_L( u \tilde u^{-1} a_i)-\ell_L(\sigma)-\ell_L(u \tilde
            u^{-1}) \\
        &= \ell_L( u \tilde u^{-1} a_i)-\ell_L(u \tilde u^{-1}) \\
        &= \ell_A( u \tilde u^{-1} a_i)-\ell_A(u \tilde u^{-1}),
\end{align*}
whence $u\tilde u^{-1}=1$, i.e. $\sigma=\tilde \sigma$.

\end{proof}

Let
$
    T = \{ \sigma \in L_{n+1} \st \des_A(\sigma) = 0 \}$.

\begin{cor}\label{COR:LnDecomposition}
1. For every $B\subseteq [n+1]$ there exists a unique $\sigma \in
T$ such that $ B=\Neg(\sigma^{-1})$.

2. For every $B\subseteq [n+1]$,
\begin{equation}
    \{\pi \in L_n \st \Neg(\pi^{-1})\subseteq B\} = \biguplus_{u \in A_{n+1}} \{ \sigma u \st \sigma \in T, \Neg(\sigma^{-1})\subseteq B
    \}   ,
\end{equation}
where $\uplus$ denotes disjoint union.
\end{cor}

\begin{cor}
Let $\pi \in L_{n+1}$, and write $\pi=\sigma u$ with $\sigma$ and
$u$ like in Lemma \ref{LE:1}. Then $\ell_L(\pi) =
\ell_A(u)+\sum_{i \in \Neg(\pi^{-1}) } i$ .
\end{cor}
\begin{proof}
By \eqref{EQ:ellL}, Lemma \ref{LE:1} and Proposition
\ref{PR:RR4.4},
\[
\begin{split}
    \ell_L(\pi)
        &= \inv(\pi) -\del_B(\pi) + \sum_{i \in \Neg(\pi^{-1}) } i\\
        &= \inv(u)   -\del_S(u) + \sum_{i \in \Neg(\pi^{-1}) } i\\
        &= \ell_A(u)  + \sum_{i \in \Neg(\pi^{-1}) } i  .
\end{split}
\]
\end{proof}

\begin{proof}[Proof of Proposition \ref{PR:1}]
From Corollary \ref{COR:LnDecomposition}, Lemma \ref{LE:1} and
Theorem \ref{TH:RR6.1},
\begin{align*}
    \sum_{\substack{ \pi \in L_{n+1} \\ \Neg(\pi^{-1})\subseteq B}} q^{\nrmaj_{L_{n+1}}(\pi)}
        &= \sum_{\substack{\sigma \in T \\ \Neg(\sigma^{-1})\subseteq B}} \sum_{u \in A_{n+1}} q^{\nrmaj_{L_{n+1}}(\sigma
            u)}  \\
        &= \sum_{\substack{\sigma \in T \\ \Neg(\sigma^{-1})\subseteq B}} \sum_{u \in A_{n+1}} q^{\rmaj_L(\sigma
            u) + \sum_{i \in \Neg((\sigma u)^{-1})} i } \\
        &= \sum_{\substack{\sigma \in T \\ \Neg(\sigma^{-1})\subseteq B}} q^{ \sum_{i \in \Neg(\sigma^{-1})} i} \sum_{u \in A_{n+1}} q^{\rmaj_A(u)} \\
        &= \sum_{ C \subseteq B } q^{ \sum_{i \in C} i} \sum_{u \in A_{n+1}} q^{\rmaj_A(u)} \\
        &= \prod_{i \in B}(1+q^i) \prod_{i=1}^{n-1}(1+q+\dots+q^{i-1}+2q^i).
\end{align*}
By similar considerations, this time invoking the other equality
in Theorem \ref{TH:RR6.1},
\begin{align*}
    \sum_{\substack{ \pi \in L_{n+1} \\ \Neg(\pi^{-1})\subseteq B}} q^{\ell_L(\pi)}
        &= \sum_{\substack{\sigma \in T \\ \Neg(\sigma^{-1})\subseteq B}} \sum_{u \in A_{n+1}}
                                                            q^{\ell_L(\sigma u)} \\
        &= \sum_{\substack{\sigma \in T \\ \Neg(\sigma^{-1})\subseteq B}} \sum_{u \in A_{n+1}}
                q^{\inv(\sigma u) + \sum_{i \in \Neg((\sigma u)^{-1})}i -\del_B(\sigma u)} \\
        &= \sum_{\substack{\sigma \in T \\ \Neg(\sigma^{-1})\subseteq B}} q^{ \sum_{i \in \Neg(\sigma^{-1})}i}
                \sum_{u \in A_{n+1}} q^{\inv(u)-\del_S(u)} \\
        &= \sum_{ C \subseteq B } q^{ \sum_{i \in C} i}
                \sum_{u \in A_{n+1}} q^{\ell_A(u)} \\
        &= \prod_{i \in B}(1+q^i) \prod_{i=1}^{n-1}(1+q+\dots+q^{i-1}+2q^i).
\end{align*}
\end{proof}

\section{Even-signed Even Permutations}\label{SEC:evensigned}
We denote by $D_n$ the group of {\it even-signed permutations},
that is the subgroup of $B_n$ consisting of all the signed
permutations having an even number of negative entries in their
window notation. Equivalently,
\[
    D_n = \{ \sigma \in B_n \st \text{$\#\Neg(\sigma^{-1})$ is even} \}.
\]

$D_n$ is a Coxeter group of type $D$, generated by $\tilde s_0,
s_1,\dots,s_{n-1}$, where $
    \tilde s_0 = [-2,-1,3,\dots,n] = s_0 s_1 s_0$.

Following Biagioli \cite{B}, we define the {\it $D$-length}  of
$\sigma \in D_n$ by
\[ \ell_D(\sigma) = \ell_B(\sigma)-\#\Neg(\sigma), \] which is
also the length of a reduced expression for $\sigma$ in the above
generators, and we let \[\dmaj(\sigma) = \maj_B(\sigma) -
\#\Neg(\sigma) + \sum_{i\in \Neg(\sigma^{-1})}i .\]

Biagioli proved the following $D_n$-analogue of MacMahon's
theorem.
\begin{prop}[See \cite{B}, Proposition 3.1]\label{PR:B3.1}
\[
    \sum_{\sigma \in D_n} q^{\dmaj(\sigma)}
        = \sum_{\sigma \in D_n} q^{\ell_D(\sigma)} .
\]
\end{prop}

Let
\[
    \drmaj_n(\sigma) = \rmaj_{B_n}(\sigma) - \#\Neg(\sigma)
    + \sum_{i \in \Neg(\sigma^{-1})}i .
\]

Since the involution $\phi$ from Lemma \ref{LE:phi} satisfies
\eqref{EQ:NegUnderPhi}, $\dmaj$ and $\drmaj_n$ are equidistributed
on $D_n$, hence we can replace $\dmaj$ with $\drmaj_n$ in
Proposition \ref{PR:B3.1}.

Let $(L \cap D)_{n+1} = L_{n+1}\cap D_{n+1}$, the group of
even-signed even permutations on $\pm 1,\dots,\pm(n+1)$, and let
\[
    \ell_{(L\cap D)}(\pi) = \ell_D(\pi) - \del_B(\pi)
\]
and
\[
    \drmaj_{(L\cap D)_{n+1}}(\pi) =
    \rmaj_{L_{n+1}}(\pi)-\#\Neg(\pi)+\sum_{i \in \Neg(\pi^{-1})}
    i.
\]

\begin{prop}\label{PR:2}
\[
    \sum_{\pi \in (L\cap D)_{n+1}} q^{\drmaj_{(L\cap
    D)_{n+1}}(\pi)} = \sum_{\pi \in (L\cap D)_{n+1} } q^{\ell_{(L\cap
        D)}(\pi)}.
\]
\end{prop}
\begin{proof}
From the definitions and from Corollary \ref{COR:inex} we have for
every $i$
\[\begin{split}
    \sum_{\substack{ \pi \in L_{n+1} \\ \#\Neg(\pi^{-1})=2i}} q^{\drmaj_{(L\cap
    D)_{n+1}}(\pi)}
        &= \sum_{\substack{ \pi \in L_{n+1} \\ \#\Neg(\pi^{-1})=2i}}
        q^{\nrmaj_{L_{n+1}}(\pi)-\#\Neg(\pi)} \\
        &= q^{-2i} \sum_{\substack{B\subseteq [n+1]\\\abs{B}=2i}} \sum_{\substack{ \pi \in L_{n+1} \\ \Neg(\pi^{-1})=B}}
        q^{\nrmaj_{L_{n+1}}(\pi)} \\
        &= q^{-2i} \sum_{\substack{B\subseteq [n+1]\\\abs{B}=2i}} \sum_{\substack{ \pi \in L_{n+1} \\ \Neg(\pi^{-1})=B}}
        q^{\ell_L(\pi)} \\
        &= \sum_{\substack{ \pi \in L_{n+1} \\ \#\Neg(\pi^{-1})=2i}}
        q^{\ell_L(\pi) -\#\Neg(\pi)} \\
        &= \sum_{\substack{ \pi \in L_{n+1} \\ \#\Neg(\pi^{-1})=2i}} q^{\ell_{(L\cap
        D)}(\pi)}.
\end{split}\]

Taking the sum over all $i$ we get the desired equality.
\end{proof}

\section*{Acknowledgements}
I would like to thank my advisor, Amitai Regev, for suggesting the topic and
for his helpful remarks on preliminary versions of this paper.

\end{document}